\documentclass[11pt]{amsart}%
\usepackage[centertags]{amsmath}
\usepackage{amsfonts}
\usepackage{amssymb}
\usepackage{amsthm}
\usepackage[english]{babel}
\usepackage[active]{srcltx}
\usepackage[colorlinks]{hyperref}
\usepackage{enumerate}
\usepackage{mathrsfs}
\usepackage{esint}
\usepackage{hyperref}
\usepackage{graphicx}%
\providecommand{\U}[1]{\protect\rule{.1in}{.1in}}
\numberwithin{equation}{section}
\newtheorem{thm}{Theorem}[section]

\newtheorem{question}[thm]{Question}

\theoremstyle{remark}
\newtheorem{rem}[thm]{Remark}

\def\R{\mathbb{R}}

\def\N{\mathbb{N}}

\begin{document}
\title[ ]{Sparse Brudnyi and John--Nirenberg Spaces}
\author{\'Oscar Dom\'inguez}
\address{O. Dom\'inguez, Departamento de An\'alisis Matem\'atico y Matem\'atica
Aplicada, Facultad de Matem\'aticas, Universidad Complutense de Madrid\\
Plaza de Ciencias 3, 28040 Madrid, Spain.}
\email{oscar.dominguez@ucm.es}
\author{Mario Milman}
\address{M. Milman, Instituto Argentino de Matematica\\
Buenos Aires\\
Argentina}
\email{mario.milman@icloud.com}
\urladdr{https://sites.google.com/site/mariomilman/}
\thanks{The first named author is supported in part by MTM2017-84058-P (AEI/FEDER, UE)}

\begin{abstract}
A generalization of the theory of Y. Brudnyi \cite{yuri}, and A. and Y.
Brudnyi \cite{BB20a}, \cite{BB20b}, is presented. Our construction connects Brudnyi's theory, which relies on local polynomial approximation,
with new results on sparse domination. In particular, we find an analogue of the
 maximal theorem for the fractional maximal function,  
solving a problem proposed by Kruglyak--Kuznetsov. Our spaces shed light
on the structure of the John--Nirenberg spaces. We show that $SJN_{p}$
(sparse John--Nirenberg space) coincides with $L^{p},1<p<\infty.$ This
characterization yields the John--Nirenberg inequality by extrapolation and is useful in the theory of commutators.

\end{abstract}
\maketitle

\section{Preamble}

The function spaces we use in Analysis can be described and characterized in
different qualitative and quantitative ways (e.g. by duality, as being part of
an interpolation scale, through the boundedness of suitable functionals, by
the rate of approximation of their elements with respect to a fixed class of
approximants, etc.). To have a complete catalog of different characterizations
at hand is of fundamental importance to understand the structure of the spaces
and facilitates their use in applications (cf. \cite{pee}, \cite{tri}).

Yuri and Alexander Brudnyi (cf. \cite{yuri, BB20a, BB20b}) have
proposed the concept of best local polynomial approximation as a unifying
characteristic to understand the structure of classical function spaces as
diverse as, $\text{BMO},$ John--Nirenberg spaces $JN_{p},$ Sobolev spaces, Besov
spaces, Morrey spaces, Jordan--Wiener spaces, etc. Their massive theory can be
seen as a complement of the theories of function spaces that have evolved
through the work of many authors, including names such as Coifman--Meyer,
Frazier--Jawerth, Peetre, Triebel (cf. \cite{fraz}, \cite{pee}, \cite{tri} and
the references therein), where the underlying unifying themes and tools are
wavelet approximations, representation theorems, maximal inequalities,
interpolation, etc. A distinguished feature of Brudnyi's constructions is the
fact that instead of explicitly defining oscillations, these appear in
Brudnyi's theory as the solution of variational problems. For example, best
local approximation by constants in $L^{1},$ amounts to replace\footnote{where
as usual, $f_{Q}$ denotes the integral average of $f$ over $Q$, i.e.,
$f_{Q}=\frac{1}{|Q|}\int_{Q}f$.} $\int_{Q}\left\vert f(x)-f_{Q}\right\vert dx$
by inf$_{c}\int_{Q}\left\vert f(x)-c\right\vert dx,$ for each cube $Q.$ Higher
order oscillations can be handled replacing constants by polynomials of a
given order and one can go on to accommodate different geometries, different
approximants, etc.

In our work on fractional maximal operators and commutators we found it very
useful to combine the constructions of Brudnyi's theory of function spaces
with new developments in Harmonic Analysis connected with Covering Lemmas,
more specifically \textquotedblleft sparseness\textquotedblright\ (cf. Lerner \cite{ler}). Our analysis led us to modify the original constructions and
formally introduce a new class of spaces that we call \emph{Sparse Brudnyi spaces}
(SB spaces). As we shall see, SB spaces can be used to provide new
characterizations of classical function spaces as well as clarify, simplify
and solve some open questions. In a different direction, our spaces build a
new bridge that should also benefit local polynomial approximation theory.

The problem of understanding the structure of the John--Nirenberg $JN_{p}$
spaces continues to attract attention to this day (cf. \cite{DHKY} for a
recent account). It was therefore surprising for us to discover that the
$SJN_{p}$ spaces (Sparse John--Nirenberg spaces) admit a simple
characterization (cf. Theorem \ref{TheoremSparseJN} below). Let $Q_{0}$ be a cube in
$\R^{n},$ then
\[
SJN_{p}(Q_{0})=\left\{
\begin{array}
[c]{lcl}%
L(\log L)(Q_{0}) & , & p=1,\\
&  & \\
L^{p}(Q_{0}) & , & 1<p<\infty,\\
&  & \\
\text{BMO}(Q_{0}) & , & p=\infty.
\end{array}
\right.
\]
Moreover, this characterization is useful in applications to a number of
problems in Analysis, where $JN_{p}$ type conditions appear naturally but, for
which, the usual embedding into weak-$L^p$, i.e.,  
\[
JN_{p} (Q_0)\subset L(p,\infty)(Q_0),\qquad 1<p<\infty,
\]
only leads to weaker inequalities (pun intended.)

As an example, we answer a question proposed by Kruglyak--Kuznetsov \cite{KK06}
concerning the fractional maximal operator. Let $M$ be the classical maximal operator, a version of the Hardy--Littlewood theorem can be
formulated as\footnote{Here the symbol $f \approx g$ indicates the existence of a universal constant $c > 0$ (independent of all parameters involved) such that $(1/c) f \leq g \leq c f$. Likewise the symbol $f \lesssim g$ will mean that there exists a universal constant $c > 0$ (independent of all parameters involved) such that $f \leq c g$.}
\begin{equation}\label{HLMaxThem}
\Vert Mf\Vert_{L^{p}(Q_{0})}\approx\Vert f\Vert_{L^{p}(Q_{0})},\qquad 1<p\leq\infty.
\end{equation}
Kruglyak--Kuznetsov \cite{KK06}, ask for an analogue for maximal
fractional operators (cf. (\ref{vea}) below). Using the parameters of
Sobolev's inequality, namely, $\lambda\in(0,n),p\in(1,\frac{n}{\lambda})$ and
$\frac{1}{q}=\frac{1}{p}-\frac{\lambda}{n},$ we only have the one direction
inequality
\begin{equation}
\Vert M_{\lambda}f\Vert_{L^{q}(\mathbb{R}^{n})}\lesssim\Vert f\Vert_{L^{p}%
(\mathbb{R}^{n})}.\label{inanalogy}%
\end{equation}
The problem thus is to find spaces that could turn (\ref{inanalogy}) into an
equivalence, and Kruglyak--Kuznetsov proposed certain spaces defined in terms of capacities. Using
SB spaces instead, we can now provide a complete answer (cf. Theorem
\ref{TheoremSparseBrudnyiV1} below)
\begin{equation*}
\Vert f\Vert_{SV_{p,q}^{k,\lambda}(Q_{0})}\approx\Vert M_{q,\lambda,Q_{0}%
}(f-P^k_{Q_0} f)\Vert_{L^{p}(Q_{0})} 
\end{equation*}
provided that $k \in \N, \lambda \in [0, n)$ and $p, q \in [1, \infty)$.


The proofs of the results announced in this note  are given in \cite{DM}.

\section{John--Nirenberg Spaces: A primer}

Initially introduced by John--Nirenberg \cite{JN61} in 1961 to study some
problems in elasticity, $\text{BMO}$ (cf. (\ref{bmo})) has been often at the center
of important developments in Analysis ever since.

It was quickly realized, probably starting with Moser \cite{mo}, that BMO
could be used to substitute $L^{\infty},$ as an endpoint space for the scale
of $L^{p}$ spaces. Indeed, many important inequalities in PDEs and Harmonic
Analysis that fail for the smaller space $L^{\infty}$ are valid for BMO.
Interpolation theorems by Stampachia \cite{stamp1, stamp2}, Campanato
\cite{campanato} and others, soon confirmed the r\^{o}le of BMO as a natural
endpoint space for the $L^{p}$ scale. A major inflection point establishing
the importance of BMO in Analysis, was Fefferman's celebrated duality
theorem \cite{fef} and the real variable approach to the theory of $H^{p}$
spaces by Fefferman and Stein \cite{fefstein}. These developments were
naturally connected with probabilistic interpretations of BMO using
martingale theory (cf. \cite{gar}). In the last 50 years or so the
applications have multiplied in many directions including the study of PDEs,
Singular Integrals (T1 theorem, Commutators), Sobolev spaces, Noncommutative
Analysis and Operator Theory. BMO conditions have been studied in many
unexpected settings, for example, in the Theory of Semigroups (probably starting
in \cite{var}) just to name a few of these developments.

An early fundamental result of BMO theory is the John--Nirenberg inequality
\cite{JN61} that gives the exponential decay of the distribution function of
elements of BMO, and implies the embedding
\begin{equation}
\text{BMO}\subset e^{L}.\label{expo}%
\end{equation}

Moreover, in the same paper, John--Nirenberg introduced the $JN_{p}$ spaces,
whose definition we now recall. Let $Q_{0}$ be a fixed cube in $\R^n$ with sides parallel to the coordinate axes, the
$JN_{p}(Q_{0})$ spaces can be defined using packings\footnote{Packings
$\Pi(Q_{0})$ are simply countable collections of subcubes of $Q_{0}$ with
pairwise disjoint interiors.} $\Pi(Q_{0}),$ as follows. For each $\pi
={\{Q_{i}\}_{i\in I}}\in\Pi(Q_{0}),$ let $f_{\pi}(x)=%
{\displaystyle\sum\limits_{i \in I}}
\Big(  \frac{1}{\left\vert Q_{i}\right\vert }\int_{Q_{i}}\left\vert
f-f_{Q_{i}}\right\vert\Big)  \mathbf{1}_{Q_{i}}(x),$ then 
$JN_{p}(Q_{0}),1\leq p<\infty,$ is defined requiring the following
functionals to be finite,
\begin{align}
\left\Vert f\right\Vert _{JN_{p}(Q_{0})} &  =\sup_{\pi\in\Pi(Q_{0})}\left\Vert
f_{\pi}\right\Vert _{L^{p}(Q_{0})}\nonumber\\
&  =\sup_{\pi\in\Pi(Q_{0})}\bigg(\sum_{i\in I}\bigg(\frac{1}{|Q_{i}|}%
\int_{Q_{i}}|f(x)-f_{Q_{i}}| \, dx\bigg)^{p}\,|Q_{i}|\bigg)^{1/p}.\label{jnpp}%
\end{align}
For $p=\infty$ we simply have,%
\begin{equation}
\left\Vert f\right\Vert _{JN_{\infty}(Q_{0})}   =\sup_{Q\subset Q_{0}}%
\frac{1}{\left\vert Q\right\vert }\int_{Q}\left\vert f(x)-f_{Q}\right\vert
dx  =\left\Vert f\right\Vert _{\text{BMO}(Q_{0})}.\label{bmo}%
\end{equation}
The decay of the distribution functions of elements of $JN_{p}(Q_0)$ can be
obtained from the embeddings%
\begin{equation}
L_{p}(Q_{0})\subset JN_{p}(Q_{0})\subset L(p,\infty)(Q_{0}),\qquad 1<p<\infty
.\label{jnp}%
\end{equation}
Bennett--DeVore--Sharpley \cite{bds}, showed that the sharp limiting result
obtains replacing $e^{L}$ by $L(\infty,\infty)=$ the rearrangement invariant
hull of BMO\footnote{The set $L(\infty,\infty)$ can be described by the
finitiness of the nonlinear functional%
\[
\left\Vert f\right\Vert _{L(\infty,\infty)}=\sup_{t > 0} \, \{f^{\ast\ast}(t)-f^{\ast
}(t)\}.
\]
} in (\ref{expo})$.$

The initial motivation for our investigation was the method to prove
(\ref{jnp}) given in \cite{GR74}, \cite{M16}, and the theory of
Garsia--Rodemich spaces\footnote{Initially (cf. \cite{M16}) these spaces were
denoted by $GaRo_{p},$ but after the definition was extended to r.i. spaces
the notation $GaRo_{X}$ was adopted (cf. \cite{astash}). In particular, in
this notation, $GaRo_{L(p,\infty)}:=GaRo_{p}$.} that has evolved since
(cf. \cite{astash} and the references therein) suggesting that perhaps there
was a simpler structure behind the $JN_{p}$ spaces.

Recall that the space $GaRo_{L(p,\infty)}(Q_{0}),\,1<p\leq\infty,$ is defined
in terms of the functional
\[
\Vert f\Vert_{GaRo_{L(p,\infty)}(Q_{0})}=\sup_{(Q_{i})_{i\in I}\in\Pi(Q_{0}%
)}\frac{\sum_{i\in I}\int_{Q_{i}}\left\vert f(x)-f_{Q_{i}}\right\vert
dx}{\big(\sum_{i\in I}|Q_{i}|\big)^{1/p^{\prime}}}.
\]
As usual, $p'$ denotes the dual exponent of $p$ given by $\frac{1}{p} + \frac{1}{p'} = 1$. 
Therefore, since for any $(Q_{i})_{i\in I}\in\Pi(Q_{0}),$ 
\begin{align*}
\frac{\sum_{i\in I}\int_{Q_{i}%
}\left\vert f(x)-f_{Q_{i}}\right\vert dx}{(\sum_{i\in I}|Q_{i}|)^{1/p^{\prime}}} &  =\frac{\sum_{i\in I}\left(  \frac{1}{\left\vert
Q_{i}\right\vert }\int_{Q_{i}}\left\vert f(x)-f_{Q_{i}}\right\vert dx\right)
\left\vert Q_{i}\right\vert ^{1/p}\left\vert Q_{i}\right\vert ^{1/p^{\prime}%
}}{(\sum_{i\in I}%
|Q_{i}|)^{1/p^{\prime}}}\\
&  \leq\left(  \sum_{i\in I}\left(  \frac{1}{\left\vert Q_{i}\right\vert }%
\int_{Q_{i}}\left\vert f(x)-f_{Q_{i}}\right\vert dx\right)  ^{p}\left\vert
Q_{i}\right\vert \right)  ^{1/p},
\end{align*}
it follows that
\begin{equation}\label{JNGR}
\left\Vert f\right\Vert _{GaRo_{L(p,\infty)}(Q_{0})}\leq\left\Vert
f\right\Vert _{JN_{p}(Q_{0})}.
\end{equation}
The remarkable result here (cf. \cite{GR74}, \cite{M16}) is that%
\begin{equation}
GaRo_{L(p,\infty)}(Q_{0})=L(p,\infty)(Q_{0}),\qquad 1<p<\infty,\label{norma-1}%
\end{equation}
and with proper definitions (cf. \cite{astash})%
\begin{equation}
GaRo_{L^{p}}(Q_{0})=\left\{
\begin{array}
[c]{cc}%
L^{p}(Q_{0}), & 1<p<\infty\\
& \\
\text{BMO}(Q_{0}), & p=\infty.
\end{array}
\right.  \label{norma-2}%
\end{equation}
Therefore the second embedding in \eqref{jnp} can be alternatively achieved as a combination of \eqref{JNGR} and \eqref{norma-1}. 

Another relevant result for us, is the classical theorem of Riesz \cite{R10}
 asserting\footnote{One estimate follows by Jensen's inequality and the other
using dyadic partitions and Lebesgue differentiation theorem.} that
\begin{equation*}
\Vert f\Vert_{L^{p}(Q_{0})}\approx\sup_{(Q_{i})_{i\in I}\in\Pi(Q_{0}%
)}\bigg(\sum_{i\in I}\bigg(\frac{1}{|Q_{i}|}\int_{Q_{i}}|f(x)|\,dx\bigg)^{p}%
\,|Q_{i}|\bigg)^{1/p}.
\end{equation*}
The previous discussion shows that this equivalence fails when dealing with
oscillations (cf. \eqref{jnpp})
\begin{equation*}
\Vert f-f_{Q_{0}}\Vert_{L^{p}(Q_{0})}\not \approx \Vert f\Vert_{JN_{p}(Q_{0}%
)}=\sup_{(Q_{i})_{i\in I}\in\Pi(Q_{0})}\bigg(\sum_{i\in I}\bigg(\frac
{1}{|Q_{i}|}\int_{Q_{i}}|f(x)-f_{Q_{i}}|\,dx\bigg)^{p}\,|Q_{i}|\bigg)^{1/p}%
.
\end{equation*}

This sets up the stage for the new ingredient in our construction.

\section{Sparse John--Nirenberg spaces}

Let $\mathcal{D}(Q_{0})$ be the
collection of all dyadic subcubes in $Q_{0}$. We say that $\mathcal{S}%
(Q_{0})\subset\mathcal{D}(Q_{0})$ is \emph{sparse} if for every $Q\in
\mathcal{S}(Q_{0})$,
\begin{equation}
\sum_{Q^{\prime}\in\mathbf{Ch}_{\mathcal{S}(Q_{0})}(Q)}|Q^{\prime}|\leq
\frac{1}{2}|Q|,\label{SparseDef}
\end{equation}
where $\mathbf{Ch}_{\mathcal{S}(Q_{0})}(Q)$ denotes the set of maximal (with
respect to inclusion) cubes in $\mathcal{S}(Q_{0})$, which are strictly
contained in $Q$.

The concept of sparse family has its roots in the classical Calder\'on--Zygmund decomposition lemma \cite[Lemma 1]{CZ}. The related sparse domination principle, which essentially establishes pointwise bounds of general Calder\'on--Zygmund operators by a supremum of a special collection of dyadic and positive operators (the so-called \emph{sparse operators}), has been recently developed into a powerful tool by Lerner \cite{ler}. Over the last few years, sparse domination has been further extended and refined to deal with many classical operators in Analysis. In this regard, we only mention \cite{lerna} and the references within. 

Let $p\in\lbrack1,\infty]$, the \emph{sparse John--Nirenberg space}
$SJN_{p}(Q_{0})$ is defined as the set of all $f\in L^{1}(Q_{0})$
such that
\begin{equation}
\Vert f\Vert_{SJN_{p}(Q_{0})}=\sup_{(Q_{i})_{i\in I}\in
\mathcal{S}(Q_{0})}\left\Vert \sum_{i\in I}\left(  \frac{1}{|Q_{i}|}%
\int_{Q_{i}}|f(x)-f_{Q_{i}}|\,dx\right)  \mathbf{1}_{Q_{i}\backslash\cup_{Q^{\prime}\in\mathbf{Ch}_{\mathcal{S}
(Q_0)}(Q_i)} Q^{\prime}}\right\Vert
_{L^{p}(Q_{0})}<\infty.\label{DefSJN}%
\end{equation}
Clearly, if $p=\infty$ then $SJN_{\infty}(Q_{0})=\text{BMO}(Q_{0})$.
Since every sparse family $(Q_i)_{i \in I} \in \mathcal{S}(Q_0)$ is, in particular, \emph{weakly sparse}, i.e., 
\begin{enumerate}[\upshape(i)]
\item the sets $E_{Q_i} = Q_{i}\backslash\cup_{Q^{\prime}\in\mathbf{Ch}_{\mathcal{S}
(Q_0)}(Q_i)} Q^{\prime}$ are pairwise disjoint,
\item $\frac{1}{2}|Q_{i}|\leq |E_{Q_i} |\leq|Q_{i}|$ (cf. \eqref{SparseDef}),
\end{enumerate}
thus we readily see that if $p\in\lbrack1,\infty)$, the expression for $\Vert
f\Vert_{SJN_{p}(Q_{0})},$ simplifies to
\begin{equation}
\Vert f\Vert_{SJN_{p}(Q_{0})}\approx\sup_{(Q_{i})_{i\in I}%
\in\mathcal{S}(Q_{0})}\left\{  \sum_{i\in I}\left(  \frac{1}{|Q_{i}|}%
\int_{Q_{i}}|f(x)-f_{Q_{i}}|\,dx\right)  ^{p}\,|Q_{i}|\right\} ^{1/p}%
.\label{SparseJNQ}%
\end{equation}
Comparing (\ref{SparseJNQ}) with (\ref{jnpp}) we see that the only difference
is that the class of admissible families of cubes is more restrictive for the
former, and therefore,
\begin{equation*}
\Vert f\Vert_{JN_{p}(Q_{0})}\leq\Vert f\Vert_{SJN_{p}(Q_{0})}.
\end{equation*}
Furthermore, $JN_{p}(Q_{0})$ is not rearrangement invariant (cf. \cite{DHKY}).
However, the situation is dramatically different for the $SJN_{p}(Q_{0})$
spaces. Indeed, we have

\begin{thm}
\label{TheoremSparseJN} Let $p \in[1, \infty]$. Then
\[
SJN_{p}(Q_{0}) = \left\{
\begin{array}
[c]{lcl}%
L^{p} (Q_{0}) & , & 1 < p < \infty,\\
&  & \\
\emph{BMO}(Q_{0}) & , & p=\infty,\\
&  & \\
L (\log L) (Q_{0}) & , & p=1.
\end{array}
\right.
\]

\end{thm}

\begin{rem}
	When we finished the first version of this paper we came across the interesting preprint \cite{ahlmo} by Airta--Hyt\"onen--Li--Martikainen--Oikari concerning mapping properties of bi-commutators on mixed norm spaces. Sparseness is also used in their method to control expressions involving oscillations. More precisely, Proposition 3.2 in \cite{ahlmo} is very close in spirit to Theorem \ref{TheoremSparseJN} with $p \in (1, \infty)$, in as much as the use of weak sparse families to control $JN_p$ conditions. But the authors use a linearization argument combined with duality, that moreover involves the Hardy--Littlewood maximal theorem. In particular, these issues preclude them from obtaining endpoint results. At any rate they do not define $SJN_p$ or, more generally, sparse variants of the Brudnyi's constructions as will appear in Section \ref{Section:Brudnyi} below. In this regard, it is important to point out that the concept of weak sparseness is not well adapted, but one needs the stronger version of sparseness provided in \eqref{SparseDef}\footnote{One only has that every weakly sparse collection decomposes into a disjoint union of finitely many sparse subcollections; cf. \cite[Lemma 6.6]{lerna}.}. 
\end{rem}

\begin{rem}
\label{RemarkBFS} The definition of the $SJN_{p}$ provided by \eqref{DefSJN} can be extended by
means of replacing $L^{p}$ with a general Banach function space (e.g.,
Lorentz space, Orlicz space, weighted Lebesgue spaces, etc.). This is not clear at all if we use instead the right-hand side of \eqref{SparseJNQ} to define $SJN_p$.
\end{rem}

As a by-product of the previous theorem and (\ref{norma-1}), (\ref{norma-2}),
we can write (\ref{jnp}) as follows
\[
GaRo_{L^{p}}(Q_{0})=SJN_{p}(Q_{0})\subset JN_{p}(Q_{0})\subset
GaRo_{L(p,\infty)}(Q_{0}),\qquad p\in(1,\infty).
\]

Theorem \ref{TheoremSparseJN} is actually a special case of a more general
result involving SB spaces, which we now introduce. 

\section{Sparse Brudnyi Spaces}\label{Section:Brudnyi}

First we review briefly the Brudnyi spaces treated in detail in \cite{yuri, BB20a, BB20b}, and for this we need to
recall the concept of $\emph{best}$ \emph{local polynomial approximation}.

For $k\in\mathbb{N}$ and $f\in L^{q}(Q_{0}),\,1\leq q\leq\infty$, we consider
the set function
\[
E_{k}(f;Q_{0})_{q}=\inf_{m\in\mathcal{P}_{k-1}^{n}}\Vert f-m\Vert_{L^{q}%
(Q_{0})},
\]
where $\mathcal{P}_{k-1}^{n}$ is the set of all polynomials in $\mathbb{R}%
^{n}$ of degree at most $k-1$. Let $\lambda\in\mathbb{R}$, and $1\leq
p\leq\infty$. We let $V_{p,q}^{k,\lambda}(Q_{0})$ denote\footnote{We warn the
reader that we have slightly changed the notation used in \cite{BB20a,
BB20b}, more precisely, the space $V_{p,q}^{k,\lambda}(Q_{0})$ defined
here corresponds with $V_{p,q}^{k,-\frac{1}{p}-\frac{\lambda}{nq}+\frac{1}%
{q}}(Q_{0})$ in those papers. We feel this parametrization helps in providing
a more clear formulation of our results. } the set of all functions $f\in
L^{q}(Q_{0})$ such that
\begin{equation}\label{BBV}
\Vert f\Vert_{V_{p,q}^{k,\lambda}(Q_{0})}=\sup_{(Q_{i})_{i\in I}\in\Pi(Q_{0}%
)}\bigg(\sum_{i\in I}(|Q_{i}|^{\frac{\lambda}{n}-\frac{1}{q}}E_{k}%
(f;Q_{i})_{q})^{p}|Q_{i}|\bigg)^{1/p}<\infty
\end{equation}
(with the usual modification if $p=\infty$). One of the main features of the
$V_{p,q}^{k,\lambda}(Q_{0})$ scale is that it can be used to provide a unified
treatment (e.g., duality assertions and structural properties) of many
classical spaces in Analysis. Since $E_{1}(f;Q)_{q}\approx(\int_{Q}%
|f-f_{Q}|^{q})^{1/q}$, we see that the spaces $JN_{p}(Q_{0}),\,p\in
(1,\infty),$ are distinguished elements of this scale, namely,
\[
JN_{p}(Q_{0})=V_{p,q}^{1,0}(Q_{0}),\qquad q\in\lbrack1,p).
\]
The list of examples of Brudnyi spaces $V_{p,q}^{k,\lambda}(Q_{0})$ also
includes, for suitable choices of the parameters, $\text{BMO}(Q_{0}), \text{BV}(Q_{0}%
),\mathring{W}^{k,p}(Q_{0})$ and $M_{q}^{\lambda}(Q_{0})$ (Morrey spaces),
among others.

Imitating the construction of the $SJN_{p}(Q_{0})$ spaces given above (cf. \eqref{DefSJN}), we can
introduce general SB spaces as follows. Let $k\in\mathbb{N},\lambda
\in\mathbb{R}$ and $p,q\in\lbrack1,\infty]$. The space $S%
V_{p,q}^{k,\lambda}(Q_{0})$ is the set of all functions $f\in L^{q}(Q_{0})$
such that
\begin{equation}
\Vert f\Vert_{SV_{p,q}^{k,\lambda}(Q_{0})}=\sup_{(Q_{i})_{i\in I}%
\in\mathcal{S}(Q_{0})}\bigg\|\sum_{i\in I}|Q_{i}|^{\frac{\lambda}{nq}-\frac
{1}{q}}E_{k}(f;Q_{i})_{q}  \mathbf{1}_{Q_{i}\backslash\cup_{Q^{\prime}\in\mathbf{Ch}_{\mathcal{S}
(Q_0)}(Q_i)} Q^{\prime}}\bigg\|_{L^{p}(Q_{0})}<\infty
.\label{DefSV}%
\end{equation}
Representative examples are given
\begin{equation}\label{SJNV}
SV_{p,1}^{1,0}(Q_{0})=SJN_{p}(Q_{0}),\qquad p\in
\lbrack1,\infty)\qquad(\text{cf. \eqref{DefSJN}}),
\end{equation}
\begin{equation*}
 SV_{\infty,1}^{1,\lambda}(Q_{0}) = \left\{
\begin{array}
[c]{lcl}%
\text{BMO}(Q_{0}) & , & \text{if } \lambda = 0,\\
&  & \\
\mathcal{C}^{-\lambda}(Q_{0}
) & , & \text{if } \lambda \in (-1, 0),\\
&  & \\
M_1^\lambda(Q_0) & , & \text{if } \lambda \in (0, n).
\end{array}
\right.
\end{equation*}
\footnote{$\mathcal{C}^{-\lambda}(Q_{0}
)$ denotes classical H\"{o}lder--Zygmund spaces} Moreover, we obviously have
\[
SV_{p,q}^{k,\lambda}(Q_{0})\subset V_{p,q}^{k,\lambda}(Q_{0}%
)\qquad\text{and}\qquad\Vert f\Vert_{V_{p,q}^{k,\lambda}(Q_{0})}\leq\Vert
f\Vert_{SV_{p,q}^{k,\lambda}(Q_{0})}.
\]

%

As a first application of the SB spaces we indicate a solution to a problem
concerning the fractional maximal operator.

For $\lambda\in\lbrack0,n)$ and $q\in\lbrack1,\infty)$, the \emph{(dyadic)
local fractional maximal operator} $M_{q,\lambda,Q_{0}}$ is defined for $f\in
L^{q}(Q_{0}),$ by
\begin{equation}
M_{q,\lambda,Q_{0}}f(x)=\sup_{\substack{Q\ni x\\Q\in\mathcal{D}(Q_{0}%
)}}\bigg(|Q|^{\frac{\lambda}{n}-1}\int_{Q}|f(y)|^{q}\,dy\bigg)^{\frac{1}{q}%
},\qquad x\in Q_{0}.\label{vea}%
\end{equation}
In particular, if $q=1$ we simply write $M_{\lambda,Q_{0}}$. If in addition
$\lambda=0,$ then the classical maximal function $M_{Q_{0}}$ is obtained, i.e., 
\[
M_{Q_{0}}f(x)=\sup_{\substack{Q\ni x\\Q\in\mathcal{D}(Q_{0})}}\frac{1}%
{|Q|}\int_{Q}|f(y)|\,dy.
\]

Kruglyak and Kuznetsov \cite[p. 310]{KK06} posed the following

\begin{question}
\label{quKK}What is a good analogue of the Hardy--Littlewood maximal theorem \eqref{HLMaxThem}
for the fractional maximal function? That is, can one construct function
spaces $X$ and $Y$ such that $L^{p} \subset X,Y\subset
L^{q}$ and $\Vert M_{\lambda, Q_0}f\Vert_{Y}\approx\Vert f\Vert_{X}$?
\end{question}

A first attempt to this question is given in Theorem 3 of
\cite{KK06} where the authors proposed an interpolation based
approach to deal with the case $Y=L^{q}(\mathbb{R}^{n},C)$ where $C$ is a
certain fractional capacity.

It turns out that the extension of Theorem \ref{TheoremSparseJN} to function
spaces with smoothness can be used to give an answer to Question \ref{quKK}.

\begin{thm}
\label{TheoremSparseBrudnyiV1} Let $k \in\mathbb{N}, \lambda\in[0, n), p
\in[1, \infty)$ and $q \in[1, \infty)$. Then
\[
S V^{k, \lambda}_{p, q}(Q_{0}) = M_{q, \lambda, Q_{0}} L^{p}%
(Q_{0}).
\]
More precisely, if $P_{Q_{0}}^{k} f \in\mathcal{P}^{n}_{k-1}$ denotes a nearly
best polynomial approximation of $f$ in $L^{q}(Q_{0})$ (i.e., $E_{k}(f;
Q_{0})_{q} \approx(\int_{Q_{0}} |f-P_{Q_{0}}^{k} f|^{q})^{1/q}$) then
\begin{equation}
\label{SparseInequality}\|M_{q, \lambda, Q_{0}}(f-P^{k}_{Q_{0}} f)\|_{L^{p}%
(Q_{0})} \leq c_{n} \|M_{Q_{0}}\|_{L^{p^{\prime}}(Q_{0}) \to L^{p^{\prime}%
}(Q_{0})} \|f\|_{S V^{k, \lambda}_{p, q}(Q_{0})}%
\end{equation}
and
\[
\|f\|_{S V^{k, \lambda}_{p, q}(Q_{0})} \leq2 \|M_{q, \lambda, Q_{0}%
}(f-P^{k}_{Q_{0}} f)\|_{L^{p}(Q_{0})}.
\]
Here, $c_{n}$ denotes a purely dimensional constant.
\end{thm}

Theorem \ref{TheoremSparseJN} is an immediate consequence of the previous
result with $\lambda=0, k=1$ and $q=1$ and the Hardy--Littlewood maximal
theorem (cf. \eqref{HLMaxThem}).

\begin{rem}
A similar comment as in Remark \ref{RemarkBFS} also applies to the previous theorem.
\end{rem}

In fact, the estimate \eqref{SparseInequality} can be sharpened if we replace
the functional $\Vert f\Vert_{SV_{p,q}^{k,\lambda}(Q_{0})}$
appearing in the right-hand side by a weaker functional involving a fractional
variant of the sparse condition \eqref{SparseDef}. This will be explained in more detail in the next section.

\section{Fractional Capacities and Sparseness}

The following definition is motivated by fractional capacities of sets of
cubes (cf. \cite{K97}). Let $\lambda\in(0,1]$. We say that $\mathcal{S}%
^{\lambda}(Q_{0})\subset\mathcal{D}(Q_{0})$ is \emph{sparse of order $\lambda
$} if for every $Q\in\mathcal{S}^{\lambda}(Q_{0})$,
\begin{equation}
\sum_{Q^{\prime}\in\mathbf{Ch}_{\mathcal{S}^{\lambda}(Q_{0})}(Q)}\left\vert
Q^{\prime}\right\vert ^{\lambda}\leq\frac{1}{2}|Q|^{\lambda}%
.\label{SparseFractDef}%
\end{equation}
Clearly, $\mathcal{S}^{1}(Q_{0})=\mathcal{S}(Q_{0})$ and
\begin{equation}
\mathcal{S}^{\lambda_{0}}(Q_{0})\subset\mathcal{S}^{\lambda_{1}}(Q_{0}%
),\qquad\lambda_{0}<\lambda_{1}\label{SparseFractEmb}%
\end{equation}
where this inclusion must be appropriately interpreted (i.e., if a given
family of cubes satisfies the condition \eqref{SparseFractDef} with
$\lambda_{0}$ then the corresponding condition with $\lambda_{1}$ also holds).
Accordingly, we can now introduce the class of function spaces
$\widetilde{SV}_{p,q}^{k,\lambda}(Q_{0}),\,\lambda\in\lbrack0,n),$
formed by all functions $f\in L^{q}(Q_{0})$ such that
\[
\Vert f\Vert_{\widetilde{SV}_{p,q}^{k,\lambda}(Q_{0})}=\sup
_{(Q_{i})_{i\in I}\in\mathcal{S}^{1-\frac{\lambda}{n}}(Q_{0})}\bigg\|\sum
_{i\in I}|Q_{i}|^{\frac{\lambda}{nq}-\frac{1}{q}}E_{k}(f;Q_{i})_{q} \mathbf{1}_{Q_{i}\backslash\cup_{Q^{\prime}\in\mathbf{Ch}_{\mathcal{S}
(Q_0)}(Q_i)} Q^{\prime}} \bigg\|_{L^{p}(Q_{0})}<\infty.
\]
Here a new phenomenon appears, namely, the supremum runs over families of sparse cubes depending on the smoothness parameter $\lambda$. In particular, this will enable us to capture better the smoothness properties of the space and is in sharp contrast with the classical constructions \eqref{BBV} where the supremum is taken with respect to all possible families of cubes. Clearly
\[
\widetilde{SV}_{p,q}^{k,0}(Q_{0})=SV_{p,q}^{k,0}(Q_{0}),
\]
and for general $\lambda\in\lbrack0,n)$, by \eqref{SparseFractEmb},
\[
SV_{p,q}^{k,\lambda}(Q_{0})\subset\widetilde{SV}%
_{p,q}^{k,\lambda}(Q_{0})\qquad\text{and}\qquad\Vert f\Vert
_{\widetilde{SV}_{p,q}^{k,\lambda}(Q_{0})}\leq\Vert f\Vert
_{SV_{p,q}^{k,\lambda}(Q_{0})}.
\]

An improvement of \eqref{SparseInequality} (cf. \eqref{SparseInequality2}
below) is contained in the following

\begin{thm}
Let $k \in\mathbb{N}, \lambda\in[0, n), p \in[1, \infty)$ and $q \in[1,
\infty)$. Then
\[
\widetilde{S V}^{k, \lambda}_{p, q}(Q_{0}) = M_{q, \lambda, Q_{0}}
L^{p}(Q_{0}).
\]
More precisely, if $P_{Q_{0}}^{k} f \in\mathcal{P}^{n}_{k-1}$ denotes a nearly
best polynomial approximation of $f$ in $L^{q}(Q_{0})$ then
\begin{equation}
\label{SparseInequality2}\|M_{q, \lambda, Q_{0}}(f-P^{k}_{Q_{0}}
f)\|_{L^{p}(Q_{0})} \leq c_{n} \|M_{Q_{0}}\|_{L^{p^{\prime}}(Q_{0}) \to
L^{p^{\prime}}(Q_{0})} \|f\|_{\widetilde{S V}^{k, \lambda}_{p,
q}(Q_{0})}%
\end{equation}
and
\[
\|f\|_{\widetilde{S V}^{k, \lambda}_{p, q}(Q_{0})} \leq2 \|M_{q,
\lambda, Q_{0}}(f-P^{k}_{Q_{0}} f)\|_{L^{p}(Q_{0})}.
\]

\end{thm}

\section{Some applications}

In this section we present further selected applications of the spaces
$SV_{p,q}^{k,\lambda}(Q_{0})$.

\subsection{A unified theory of commutators.}

To simplify our presentation, the results given in this section are only
stated for the Hilbert transform $H$ on $\mathbb{R}$, but corresponding
results for smooth Calder\'{o}n--Zygmund operators on $\mathbb{R}^{n}$ also
hold. An important family of commutators in Complex Analysis, Nonlinear PDE's,
Operator Theory and Interpolation Theory is given by
\[
\lbrack H,b]f=H(bf)-bH(f),\,\qquad b\in L_{\text{loc}}^{1}(\mathbb{R}).
\]
The mapping properties of this operator between Lebesgue spaces are collected
in the following

\begin{thm}
\label{TheoremCommutator} Let $1 < p, q < \infty$. Then
\[
[H, b]: L^{p}(\mathbb{R}) \to L^{q}(\mathbb{R})
\]
if and only if

\begin{enumerate}
[\upshape(i)]

\item $p=q$ and $b \in\emph{BMO}(\mathbb{R})$,

\item $p < q$ and $b \in\mathcal{C}^{\alpha}(\mathbb{R})$ with $\alpha=
\frac{1}{p} - \frac{1}{q}$,

\item $p > q$ and $b \in L^{r}(\mathbb{R})$ with $\frac{1}{r} = \frac{1}%
{q}-\frac{1}{p}$.
\end{enumerate}
\end{thm}

The regime $p\leq q$ in the previous theorem is classical. Specifically, the
diagonal case $p=q$ corresponds to the celebrated Coifman--Rochberg--Weiss
theorem \cite{CR76} (with \cite{N57} as a forerunner), while the case
$p<q$ is due to Janson \cite{J78}. On the other hand, the case $q<p$ was
only achieved recently by Hyt\"{o}nen \cite{H} and it is intimately
connected with the Jacobian equation. Note that the non-trivial assertion in
this case concerns the necessity of $b\in L^{r}(\mathbb{R})$ and, in
particular, it shows that the cancellation inherited to the commutator does
not play a role. The results of our
paper can now be applied to give a unified treatment of these three cases and, in particular, they show that the appropriate function spaces in commutator theorems are
$SV_{p,q}^{k,\lambda}(\mathbb{R}^{n})$, which are naturally defined
from \eqref{DefSV}. For instance, Theorem \ref{TheoremCommutator} can now be
rewritten as follows.



\begin{thm}
\label{TheoremCRWRef} Assume $1 < p, q < \infty$. Let $\lambda= -
\big(\frac{1}{p} - \frac{1}{q} \big)_{+}$ and $\frac{1}{r} = \big(\frac{1}%
{q}-\frac{1}{p}\big)_{+}$. Then
\[
[H, b]: S V^{1, 0}_{p, 1}(\mathbb{R}) \to SV^{1, 0}_{q,
1}(\mathbb{R}) \iff b \in SV^{1, \lambda}_{r, 1}(\mathbb{R}).
\]

\end{thm}

The previous statement should be adequately interpreted since $S
V^{1, 0}_{p, 1}(\mathbb{R})$ coincides with $L^{p}(\mathbb{R})$ modulo
constants (cf. Theorem \ref{TheoremSparseJN} and \eqref{SJNV}).

The reformulation of the commutator theorem given in Theorem
\ref{TheoremCRWRef} paves the way to further lines of research (cf.
\cite{DM}). For instance, it is natural to investigate what is the role played
by the parameters $k,\lambda,p$ and $q$ in $SV_{p,q}^{k,\lambda
}(\mathbb{R}^{n})$ in the commutator theorem; the pair $(L^{p},L^{q})$ can be
replaced by more general pairs of Banach function spaces; the Hilbert
transform can be replaced by another classical operators in Analysis such as
maximal functions (cf. \cite{BMR00}.)

\subsection{John--Nirenberg inequalities.}

We can give an elementary proof of the John--Nirenberg embedding \eqref{expo}
via the spaces $SJN_{p}(Q_{0})$. Indeed, given $p>1$, it follows
from Theorem \ref{TheoremSparseJN} (cf. also Theorem
\ref{TheoremSparseBrudnyiV1} with $k=1, \lambda=0$ and $q=1$) that
\[
\Vert f-f_{Q_{0}}\Vert_{L^{p}(Q_{0})}\leq c_{n}\Vert M_{Q_{0}}\Vert
_{L^{p^{\prime}}(Q_{0})\rightarrow L^{p^{\prime}}(Q_{0})}\Vert f\Vert
_{SJN_{p}(Q_{0})}.
\]
Furthermore, a well-known interpolation argument yields that
\[
\Vert M_{Q_{0}}\Vert_{L^{p^{\prime}}(Q_{0})\rightarrow L^{p^{\prime}}(Q_{0}%
)}\leq\frac{p^{\prime}}{p^{\prime}-1}.
\]
Note that $\frac{p^{\prime}}{p^{\prime}-1}=p$. Therefore
\[
\Vert f-f_{Q_{0}}\Vert_{L^{p}(Q_{0})}\lesssim p\,\Vert f\Vert_{S%
JN_{p}(Q_{0})}\leq p|Q_{0}|^{1/p}\Vert f\Vert_{\text{BMO}(Q_{0})}%
\]
and by classical extrapolation we arrive at \eqref{expo}.

The above argument can be applied \emph{mutatis mutandis} to derive John--Nirenberg-type inequalities for Morrey spaces. Further details are left to the reader. 

\subsection{Sobolev inequalities}

Theorem \ref{TheoremSparseBrudnyiV1} can be applied to give an elementary
proof of the local counterpart of the Sobolev inequality \eqref{inanalogy}.
Indeed, let $\lambda\in(0, n), \, p \in(1, \frac{n}{\lambda})$ and $\frac
{1}{q} = \frac{1}{p} - \frac{\lambda}{n}$, we have
\begin{align*}
\|M_{\lambda, Q_{0}} (f-f_{Q_{0}})\|_{L^{q}(Q_{0})}  &  \approx
\|f\|_{S V^{1, \lambda}_{q, 1}(Q_{0})}\\
&  = \sup_{(Q_{i})_{i \in I} \in\mathcal{S}(Q_{0})} \bigg\| \sum_{i \in I}
|Q_{i}|^{\frac{\lambda}{n}-1} E_{1}(f;Q_{i})_{1}  \mathbf{1}_{Q_{i}\backslash\cup_{Q^{\prime}\in\mathbf{Ch}_{\mathcal{S}
(Q_0)}(Q_i)} Q^{\prime}}
\bigg\|_{L^{q}(Q_{0})}\\
&  \leq\sup_{(Q_{i})_{i \in I} \in\mathcal{S}(Q_{0})} \bigg( \sum_{i \in I}
(|Q_{i}|^{\frac{\lambda}{n}-1} E_{1}(f;Q_{i})_{1})^{q} |Q_{i}| \bigg)^{1/q}\\
&  \leq\sup_{(Q_{i})_{i \in I} \in\mathcal{S}(Q_{0})} \bigg( \sum_{i \in I}
(|Q_{i}|^{-1} E_{1}(f;Q_{i})_{1})^{p} |Q_{i}| \bigg)^{1/p}\\
&  = \|f\|_{SJN_{p}(Q_{0})} \approx\|f-f_{Q_{0}}\|_{L_{p}(Q_{0})}.
\end{align*}

\end{document}